\documentclass{amsproc}
\usepackage{amsfonts}
\usepackage{graphicx}
\usepackage{amsmath}
\usepackage{amssymb}
\usepackage{amsthm}
\usepackage{enumerate}
\usepackage[all, cmtip]{xy} 
\newtheorem{theorem}{Theorem}

\newtheorem{corollary}[theorem]{Corollary}

\newtheorem*{mydef}{Definition}

\newtheorem{example}[theorem]{Example}

\newtheorem{lemma}[theorem]{Lemma}

\newtheorem{proposition}[theorem]{Proposition}

\newcommand{\Ker}{\operatorname{Ker}}
\newcommand{\End}{\operatorname{End}}
\begin{document}
\title{Automorphism-invariant modules}
\author{Pedro A. Guil Asensio}
\thanks{The first author was partially supported by the DGI (MTM2010-20940-C02-02)
and by the Excellence Research Groups Program of the S\'eneca Foundation of the Region of Murcia. Part of
the sources of both institutions come from the FEDER funds of the European Union.}
\address{Departamento de Mathematicas, Universidad de Murcia, Murcia, 30100, Spain}
\email{paguil@um.es}
\author{Ashish K. Srivastava}
\address{Department of Mathematics and Computer Science, St. Louis University, St.
Louis, MO-63103, USA}
\email{asrivas3@slu.edu}
\dedicatory{Dedicated to the memory of Carl Faith}

\begin{abstract}
A module is called automorphism-invariant if it is invariant under any automorphism of its injective envelope. In this survey article we present the current state of art dealing with such class of modules.  
\end{abstract}

\maketitle

\section{Introduction.}

\noindent Johnson and Wong \cite{JW} proved that a module $M$ is invariant under any endomorphism of its injective envelope if and only if any homomorphism from a submodule of $M$ to $M$ can be extended to an endomorphism of $M$. A module satisfying any of the above mentioned equivalent conditions is called a {\it quasi-injective module}. Clearly any injective module is quasi-injective. Most of the attempts of generalizing notions of injectivity or quasi-injectivity have focussed on relaxing conditions of lifting property of homomorphisms. For example, a module $M$ was called {\it pseudo-injective} by Jain et al in \cite{JS} if every monomorphism from a submodule of $M$ to $M$ extends to an endomorphism of $M$ (see \cite{AEJ}, \cite{hai}). Dickson and Fuller were first to generalize the other aspect of quasi-injective modules that these are precisely the modules that are invariant under endomorphisms of their injective envelope. Dickson and Fuller studied modules that are invariant under automorphisms of their injective envelopes in \cite{DF} for the particular case of finite-dimensional algebras over fields $\mathbb F$ with more than two elements. But recently this notion has been studied for modules over any ring. A module $M$ which is invariant under automorphisms of its injective envelope has been called an {\it automorphism-invariant module} in \cite{LZ}. The dual notion has been defined in \cite{SS1}. 

Let $M$ be an automorphism-invariant module and $M=A\oplus B$. Let $E(A)$ and $E(B)$ be injective envelopes of $A$ and $B$, respectively. Then $E(M)=E(A)\oplus E(B)$. Let $f:E(A)\rightarrow E(A)$ be any automorphism and consider the diagonal automorphism 
$(f,1_{E(B)}):E(A)\oplus E(B)\rightarrow E(A)\oplus E(B)$. As $M$ is automorphism-invariant, we get that 
$(f,1_{E(B)})(M)\subseteq M$. But this means that $f(A)\subseteq A$ by construction. Thus it follows that 

\begin{lemma} \cite{LZ}
A direct summand of an automorphism-invariant module is automorphism-invariant.
\end{lemma}     

\begin{lemma} \cite{LZ}
If for two modules $M_1$ and $M_2$, $M_1\oplus M_2$ is automorphism-invariant, then $M_1$ is $M_2$-injective and $M_2$ is $M_1$-injective.
\end{lemma} 

\noindent In \cite{LZ} it was also noted that a module $M$ is automorphism-invariant if and only if every isomorphism between any two essential submodules of $M$ extends to an automorphism of $M$. Thus it follows that any pseudo-injective module is automorphism-invariant. Lee and Zhou in \cite{LZ} asked if automorphism-invariant modules are the same as pseudo-injective modules. This was answered in the affirmative by Er, Singh and Srivastava in \cite{ESS} where it is shown that

\begin{theorem}
For a module $M$, the following are equivalent;
\begin{enumerate}[(a)]
\item $M$ is invariant under any automorphism of its injective envelope.
\item Any monomorphism from a submodule of $M$ to $M$ extends to an endomorphism of $M$.
\end{enumerate}
\end{theorem} 

\noindent We would like to emphasize here that although automorphism-invariant modules and pseudo-injective modules are the same, the idea of looking at invariance under automorphisms of injective envelope is much simpler to work with as compared to earlier used definition of (pseudo-injective and hence) automorphism-invariant modules. Teply \cite{Teply} gave a very difficult construction for modules that are automorphism-invariant but not quasi-injective. But this new way of looking at such modules gives us examples that are much easier to construct. 

\bigskip

\noindent {\bf Example.} \cite{ESS} Let $R$ be the ring of all eventually constant sequences $(x_n)_{n \in \mathbb{N}}$ of elements in  $\mathbb{F}_2$. Then $E(R_R)=\prod_{n\in \mathbb{N}}\mathbb{F}_2$, which has only one automorphism, namely the identity automorphism. Thus, $R_R$ is automorphism-invariant but it is not quasi-injective. ~\hfill$\square$

\bigskip

\noindent As a consequence of the fact that automorphism-invariant modules and pseudo-injective modules are the same, it follows easily that 

\begin{corollary} \cite{ESS}
Any automorphism-invariant module $M$ satisfies the C2 property, that is, every submodule of $M$ isomorphic to a direct summand of $M$ is itself a direct summand of $M$. 
\end{corollary}
 
 \bigskip
 
 \section{When is an automorphism-invariant module quasi-injective}
 
\noindent As automorphism-invariant modules generalize the notion of quasi-injective modules, it is natural to ask when is an automorphism-invariant module quasi-injective. This question has a natural connection to the problem of characterizing when endomorphism ring of a module is additively generated by its unit elements. If each endomorphism of the injective envelope of a module $M$ is a sum of automorphisms then clearly $M$ is automorphism-invariant if and only if $M$ is quasi-injective. Laszlo Fuchs raised the question of determining when an endomorphism ring is generated additively by automorphisms. For abelian groups, this question has been studied by many authors including Hill \cite{H} and Stringall \cite{St}. 

The structure theory of von Neumann regular right self-injective rings yields that any von Neumann regular right self-injective ring $R$ is a direct product of an abelian regular ring and a product of proper matrix rings over elementary divisor rings. Using this observation Khurana and Srivastava proved that each element of a right self-injective ring $R$ is the sum of two units if and only if $R$ has no homomorphic image isomorphic to the field of two elements $\mathbb F_2$ (see \cite{KS1}, \cite{KS2}). Thus it follows that each endomorphism of a continuous (and hence also of injective and quasi-injective) module $M$ is the sum of two automorphisms if and only if $\End(M)$ has no homomorphic image isomorphic to $\mathbb{F}_{2}$. Guil Asensio and Herzog proved that if $M$ is a flat cotorsion (in particular, pure injective) module, then $\End(M)/J(\End(M))$ is a von Neumann regular right self-injective  ring \cite{GH}. Thus each endomorphism of a flat cotorsion module $M$ is the sum of two automorphisms if and only if $\End(M)$ has no homomorphic image isomorphic to $\mathbb{F}_{2}$. 

In \cite{GS2}, Guil Asensio and Srivastava showed that if $M$ is a right $R$-module such that $\End(M)$ has no homomorphic image isomorphic to $\mathbb{F}_2$, then $\End(E(M))$ has no homomorphic image isomorphic to $\mathbb{F}_2$ either. Thus it follows that

\begin{theorem} \cite{GS2} \label{quasi}
If $M$ is a right $R$-module such that $\End(M)$ has no homomorphic image isomorphic to $\mathbb{F}_2$, then $M$ is quasi-injective if and only $M$ is automorphism-invariant.
\end{theorem}

Using the above theorem and the observation that if $R$ is any ring and $S$, a subring of its center $Z(R)$ such that $\mathbb{F}_2$ does not admit a structure of right $S$-module, then for any right $R$-module $M$, the endomorphism ring $\End(M)$ has no homomorphic image isomorphic to $\mathbb{F}_2$, Guil Asensio and Srivastava proved in \cite{GS2} that

\begin{theorem} \cite{GS2}
Let $A$ be an algebra over a field $\mathbb F$ with more than two elements. Then any right $A$-module $M$ is automorphism-invariant if and only if $M$ is quasi-injective.
\end{theorem}

This extends a result of Dickson and Fuller \cite{DF} where they proved that if $R$ is a finite-dimensional algebra over a field $\mathbb F$ with more than two elements then an indecomposable automorphism-invariant right $R$-module must be quasi-injective. The following example is from \cite{SS1} and it shows why we need to assume that the field $\mathbb F$ has more than two elements in the above theorem. 

\bigskip

\noindent {\bf Example.} Let $R=\left[ 
\begin{array}{ccc}
\mathbb F_2 & \mathbb F_2 & \mathbb F_2 \\ 
0 & \mathbb F_2 & 0 \\ 
0 & 0 & \mathbb F_2  \\ 
\end{array}
\right] $ where $\mathbb F_2$ is the field with two elements. 

\bigskip

\noindent This ring $R$ is an example of a finite-dimensional algebra over a field with two elements such that there exists an indecomposable right $R$-module $e_{11}R$ which is automorphism-invariant but not quasi-injective. ~\hfill$\square$    

\bigskip

\noindent As a consequence of the above theorem, we also have the following

\begin{corollary} \cite{GS2}
If $A$ is an algebra over a field $\mathbb F$ with more than two elements such that $A$ is automorphism-invariant as a right $A$-module, then $A$ is right self-injective.
\end{corollary}

\noindent Thus, in particular, we have the following

\begin{corollary} \cite{GS2}
Let $K[G]$ be an automorphism-invariant group algebra, where $K$ is a field with more than two elements. Then the group $G$ must be finite.
\end{corollary}

We have already seen that an automorphism-invariant module $M$ satisfies the property C2. Now if we assume, in addition, that $M$ satisfies the property C1, then $M$ is a continuous module and hence $M$ is invariant under any idempotent endomorphism of $E(M)$. Since $E(M)$ is a clean module, being an injective module, any endomorphism of $E(M)$ is a sum of an idempotent endomorphism and an automorphism. Thus, as a CS automorphism-invariant module $M$ is invariant under both idempotent endomorphisms and automorphisms of $E(M)$, we have the following.

\begin{proposition} \cite{LZ}
If $M$ is a CS automorphism-invariant module, then $M$ is quasi-injective.
\end{proposition} 

The above result shows, in particular, that if a continuous module $M$ is not quasi-injective then $M$ cannot be automorphism-invariant. Also, it shows that if $M$ is an automorphism-invariant module which is not quasi-injective then $M$ cannot be a continuous module.
\bigskip

\section{Endomorphism rings and structure of automorphism-invariant modules}

\noindent Faith and Utumi studied endomorphism ring of a quasi-injective module and extended the result of Warfield about endomorphism ring of an injective module by proving that if $M$ is a quasi-injective module and $R=\End(M)$, then $J(R)$ consists of all endomorphisms of $M$ having essential kernel and $R/J(R)$ is a von Neumann regular ring \cite{FU}. Later, Osofsky \cite{O} proved that $R/J(R)$ is right self-injective too.

In \cite{GS1} Guil Asensio and Srivastava extended the above mentioned result of Faith and Utumi to automorphism-invariant modules by proving that 

\begin{proposition} \cite{GS1}
Let $M$ be an automorphism-invariant module and $R=\End(M)$. Then $\Delta=\{f \in R: \Ker(f)\subseteq_e M\}$ is the Jacobson radical of $R$, $R/J(R)$ is a von Neumann regular ring and idempotents lift modulo $J(R)$.
\end{proposition}

\noindent Note that in the above proposition, unlike quasi-injective modules, $R/J(R)$ is not necessarily right self-injective. In the example where $R$ is the ring of all eventually constant sequences $(x_n)_{n \in \mathbb{N}}$ of elements in  $\mathbb{F}_2$, as $R$ is von Neumann regular, $J(R)=0$ and thus $R/J(R)$ is not self-injective. It was shown in \cite[Corollary 3.8]{LZ2} that the endomorphism ring of a pseudo-injective module is von Neumann regular modulo its Jacobson radical.

In \cite{ESS} and \cite{GKS}, the structure of an automorphism-invariant module is studied and it is shown that

\begin{theorem} (\cite{ESS}, \cite{GKS}) \label{decomp} 
Let $M$ be an automorphism-invariant module. Then $M$ has a decomposition $M=A\oplus B$ where $A$ is quasi-injective and $B$ is square-free.   
\end{theorem}

\noindent Recall that a module $M$ is called {\it square-free} if $M$ does not contain a nonzero submodule $N$ isomorphic to $X\oplus X$ for some module $X$. Since $B$ is square-free, all the idempotents in $\End(B)/J(\End(B))$ are central (see \cite[Lemma 3.4]{MM}). Consequently, we have that

\begin{theorem} \cite{GKS} \label{struct}
Let $M$ be an automorphism-invariant module. Then \[\End(M)/J(\End(M)) \cong R_1 \times R_2\] where $R_1$ is a von Neumann regular right self-injective ring and $R_2$ is an abelian regular ring.  
\end{theorem} 

\noindent As a consequence of this structure of automorphism-invariant modules, it may be deduced that

\begin{theorem} \cite{ESS}
If $R$ is a prime right non-singular, right automorphism-invariant ring, then $R$ is right self-injective.
\end{theorem}

\noindent In particular, this answers the question raised by Clark and Huynh in \cite[Remark 3.4]{CH1}.

\begin{corollary} \cite{ESS}
A simple right automorphism-invariant ring is right self-injective.
\end{corollary}

\noindent To understand the square-free part of automorphism-invariant modules, note that if $M$ is a square-free automorphism-invariant module, then $\End(M)/J(\End(M))$ is an abelian regular ring. Therefore, $\End(M)/J(\End(M))$ and hence $\End(M)$ is both right as well as left quasi-duo. Recall that a ring $R$ is called a right {\it quasi-duo ring} if every maximal right ideal of $R$ is two-sided. It is not known whether quasi-duo rings are left-right symmetric. 

Now, assume $N$ is an automorphism-invariant module such that $\End(N)$ is left quasi-duo. We claim that in this situation $N$ is square-free. Assume to the contrary that $N=N_1\oplus N_2\oplus N_3$ with $N_1\cong N_2$. Let $I$ be a maximal left ideal of $\End(N)$ containing $a=i_{N_2\oplus N_3} \pi_{N_2\oplus N_3}$ where $i_{N_2\oplus N_3}:N_2\oplus N_3\rightarrow N$ and $\pi_{N_2\oplus N_3}:N\rightarrow N_2\oplus N_3$ are structural injection and projection. Since $\End(N)$ is left quasi-duo, $I$ is a two-sided ideal of $\End(N)$. Let $\varphi:N_1\rightarrow N_2$ be an isomorphism. Define $f=i_{N_2} \varphi \pi_{N_1}$ and $g=i_{N_1} \varphi^{-1} \pi_{N_2} + i_{N_3}\pi_{N_3}$. Then $af+ga\in I$ and it is a monomorphism. Since $N$ is automorphism-invariant, there exists an $h\in \End(N)$ such that $h(af+ga)=1$. Thus $I=\End(N)$, a contradiction. Hence $N$ is square-free. Thus, we have

\begin{theorem}
Let $M$ be an automorphism-invariant module. Then $M$ is square-free if and only if $\End(M)$ is left (and hence right) quasi-duo. 
\end{theorem}

\bigskip

\section{Properties of automorphism-invariant modules}

\noindent The notion of exchange property for modules was introduced by Crawley and J\'{o}nnson \cite{CJ}. A right $R$-module $M$ is said to satisfy the {\it exchange property} if for every right $R$-module $A$ and any two direct sum decompositions $A=M^{\prime}\oplus N=\oplus_{i\in \mathcal I}A_{i}$ with $M^{\prime} \simeq M$, there exist submodules $B_i$ of $A_i$ such that $A=M^{\prime} \oplus (\oplus_{i \in \mathcal I}B_i)$.
If this hold only for $|\mathcal I|<\infty$, then $M$ is said to satisfy the finite exchange property. Crawley and J\'{o}nnson raised the question whether the finite exchange property always implies the full exchange property but this question is still open. 

A ring $R$ is called an {\it exchange ring} if the module $R_R$ (or $_RR$) satisfies the (finite) exchange property. Warfield \cite{Warfield2} showed that exchange rings are left-right symmetric and that a module $M$ has the finite exchange property if and only if $\End(M)$ is an exchange ring.     

Warfield \cite{Warfield1} proved that injective modules satisfy the full exchange property. This was later extended by Fuchs \cite{Fuchs} who showed that quasi-injective modules too satisfy the full exchange property. So it is natural to ask whether automorphism-invariant modules satisfy the exchange property. We have already discussed that if $M$ is an automorphism-invariant module, then $\End(M)/J(\End(M))$ is a von Neumann regular ring and idempotents lift modulo $J(\End(M)))$. Thus by \cite[Proposition 1.6]{Nichol}, it follows that $\End(M)$ is an exchange ring. This shows that $M$ has the finite exchange property. Keeping in mind the decomposition $M=A\oplus B$ where $A$ is quasi-injective and $B$ is a square-free module and the facts that for a square-free module, the finite exchange property implies the full exchange property and that every quasi-injective module satisfies the full exchange property, we obtain

\begin{theorem} \cite{GS1}
An automorphism-invariant module satisfies the exchange property.
\end{theorem}

\noindent Nicholson introduced the notion of clean rings in \cite{Nichol}. A ring $R$ is called a {\it clean ring} if each element $a\in R$ can be expressed as $a=e+u$ where $e$ is an idempotent in $R$ and $u$ is a unit in $R$. It is not difficult to see that clean rings are exchange. A module $M$ is called a clean module if $\End(M)$ is a clean ring. For example, continuos modules are known to be clean. 

Nicholson proved that an exchange ring with central idempotents is a clean ring. Thus, in particular, an abelian regular ring is clean. It is known that any right self-injective ring is clean. Now, for any automorphism-invariant module $M$, we have that $\End(M)/J(\End(M)) \cong R_1 \times R_2$ where $R_1$ is a von Neumann regular right self-injective ring and $R_2$ is an abelian regular ring. Since both $R_1$ and $R_2$ are clean rings and direct product of clean rings is clean, it follows that $\End(M)/J(\End(M))$ is a clean ring. It is known that if $R/J(R)$ is a clean ring and idempotents lift modulo $J(R)$, then $R$ is a clean ring. Because idempotents lift modulo $J(\End(M))$, we have the following result.     

\begin{theorem} \cite{GS1}
Automorphism-invariant modules are clean.
\end{theorem}

\noindent A module $M$ is called {\it directly-finite} if $M$ is not isomorphic to a proper summand of itself. A ring $R$ is called directly-finite if $xy=1$ implies $yx=1$ for any $x, y\in R$. It is well-known that a module $M$ is directly-finite if and only if its endomorphism ring $\End(M)$ is directly-finite. 

A module $M$ is said to have the {\it cancellation property} if whenever $M\oplus A\cong M\oplus B$, then $A\cong B$. A module $M$ is said to have the {\it internal cancellation property} if whenever $M=A_1\oplus B_1\cong A_2\oplus B_2$ with $A_1\cong A_2$, then $B_1\cong B_2$. A module with the cancellation property always satisfies the internal cancellation property but the converse need not be true, in general. Fuchs \cite{Fuchs} proved that if $M$ is a module with the finite exchange property, then $M$ has the cancellation property if and only if $M$ has the internal cancellation property. 

A module $A$ is said to have the {\it substitution property} if for every module $M$ with decompositions $M=A_1\oplus H=A_2\oplus K$ with $A_1\cong A\cong A_2$, there exists a submodule $C$ of $M$ (necessarily $\cong A$) such that $M=C\oplus H=C\oplus K$. In general, we have the following relation among these notions;\\
substitution$\implies$ cancellation$\implies$internal cancellation$\implies$directly-finite.

\bigskip

\noindent In \cite{GS1} Guil Asensio and Srivastava proved that if $M$ is a directly-finite automorphism-invariant module, then $E(M)$ is also directly-finite. Using this they showed that if $M$ is a directly-finite automorphism-invariant module, then $\End(M)/J(\End(M))$ is unit-regular. Consequently, it follows that

\begin{theorem} \cite{GS1}
Let $M$ an automorphism-invariant module. Then the following are equivalent;
\begin{enumerate}[(a)]
\item $M$ is directly finite.
\item $M$ has the internal cancellation property.
\item $M$ has the cancellation property.
\item $M$ has the substitution property.
\end{enumerate}
\end{theorem}

\bigskip

\section{Automorphism-invariant Leavitt path algebras}

\bigskip

\noindent In a recent preprint \cite{ARS}, automorphism-invariant Leavitt path algebras have been studied. Let $K$ be a field and $E$ be an arbitrary directed graph. Let $E^0$ be the set of vertices, and $E^1$ be the set of edges of directed graph $E$. Consider two maps $r: E^1 \rightarrow E^0$ and $s:E^1 \rightarrow E^0$. For any edge $e$ in $E^1$, $s(e)$ is called the {\it source} of $e$ and $r(e)$ is called the {\it range} of $e$. If $e$ is an edge starting from vertex $v$ and pointing toward vertex $w$, then we imagine an edge starting from vertex $w$ and pointing toward vertex $v$ and call it the {\it ghost edge} of $e$ and denote it by $e^*$. We denote by $(E^1)^*$, the set of all ghost edges of directed graph $E$. If $v \in E^0$ does not emit any edges, i.e. $s^{-1}(v) = \emptyset$, then $v$ is called a {\it sink} and if $v$ emits an infinite number of edges, i.e. $|s^{-1}(v)| = \infty$, then $v$ is called an {\it infinite emitter}. If a vertex $v$ is neither a sink nor an infinite emitter, then $v$ is called a {\it regular vertex}. 

\bigskip

\noindent The {\it Leavitt path algebra} of $E$ with coefficients in $K$, denoted by $L_K(E)$, is the $K$-algebra generated by the sets $E^0$, $E^1$, and $(E^1)^*$, subject to the following conditions:
	\begin{enumerate}
		\item[(A1)] $v_iv_j = \delta_{ij} v_i$ for all $v_i, v_j \in E^0$.
		\item[(A2)] $s(e)e = e = er(e)$ and $r(e)e^* = e^* = e^*s(e)$ for all $e$ in $E^1$.
		\item[(CK1)] $e_i^*e_j = \delta_{ij} r(e_i)$ for all $e_i, e_j \in E^1$.
		\item[(CK2)]  If $v \in E^0$ is any regular vertex, then $v = \sum_{\{e \in E^1: s(e) = v\}} ee^*$.
	\end{enumerate}
Conditions (CK1) and (CK2) are known as the {\it Cuntz-Krieger relations}. If $E^0$ is finite, then $\sum\limits_{v_i \in E^0} v_i$ is an identity for $L_K(E)$ and if $E^0$ is infinite, then $E^0$ generates a set of local units for $L_K(E)$.	 The reader is referred to \cite{AA1}, \cite{AA2} and \cite{AAJZ} for more details on Leavitt path algebras.

\bigskip

\noindent A vertex $v\in E^0$ is called a {\it bifurcation} if $|s^{-1}(v)|\ge 2$. An {\it infinite path} $\gamma$ is a sequence of edges $e_1 e_2 \ldots e_n \ldots $ such that $r(e_i)=s(e_{i+1})$ for each $i\in \mathbb N$. The infinite path $\gamma$ is called an {\it infinite sink} if there are no bifurcations nor cycles at any vertex in the path. We say that an infinite path $p$ ends in a sink if there exists an infinite sink $\gamma$ and edges $e_1, \ldots e_n\in E^1$ such that $p=e_1e_2\ldots e_n \gamma$.

\bigskip 

\noindent  Aranda Pino et al. \cite{ARS} have recently proved the following

\begin{theorem} \cite{ARS}
Let $E$ be a graph such that no vertex in $E$ emits infinitely many edges and let $R=L_K(E)$. Then the following are equivalent;
\begin{enumerate}[(a)]
\item $R$ is automorphism-invariant as a right $R$-module.
\item $R$ is right continuous.
\item $R$ is left and right self-injective and von Neumann regular.
\item $E$ is acyclic and every infinite path ends in a sink.
\end{enumerate}
\end{theorem} 

\bigskip

\section{Modules invariant under automorphisms of their pure-injective envelope}

\bigskip

\noindent In \cite{GKS}, Guil Asensio, Keskin T\"ut\"unc\"u and Srivastava have developed  a general theory of modules which are invariant under automorphisms of their envelope. The theory applies, in particular, to the modules which are invariant under automorphisms of their pure-injective envelope. We would like to finish this paper by developing in this section a more direct  approach to the study of modules which are invariant under automorphisms of their pure injective envelope, based on the results we have included in the previous sections.

Let us recall that a short exact sequence
$$0\rightarrow N\rightarrow M\rightarrow M/N\rightarrow 0$$
in ${\rm Mod}$-$R$ is called {\em pure} if the induced sequence
$$0\rightarrow N\otimes_RX\rightarrow M\otimes_RX\rightarrow (M/N)\otimes_RX\rightarrow 0$$
remains exact in $Ab$ for any left $R$-module $X$, equivalently, if any finitely presented module $F$
is projective with respect to it. And a module  $E\in{\rm Mod}$-$R$ is called {\em
pure-injective} if it is injective with respect to any pure-exact
sequence.

When dealing with pure-injectivity, the so-called `functor category technique' is quite frequently used. Let us briefly
explain the basic ideas behind this technique. An abelian category $\mathcal C$ is called a {\em Grothendieck category} if $\mathcal C$ has coproducts, exact direct limits and a generator set of objects. And, a Grothendieck category $\mathcal C$ is called {\em
locally finitely presented} if it has a generator set $\{C_i\}_{i\in \mathcal I}$ consisting of 
finitely presented objects. Recall that an object $C\in \mathcal C$ is called finitely presented if the functor ${\rm Hom}_{\mathcal C}(C,-):\mathcal{C}\rightarrow Ab$ 
commutes with direct limits. Every locally finitely presented Grothendieck category $\mathcal C$ has enough injective objects
and every object $C\in \mathcal C$ can be essentially embedded in
an injective object $E(C)$, called the injective envelope of $C$
(see e.g. \cite{Si}).

It is well-known (see e.g. \cite{CB, Si}) that one can associate to any module category Mod-$R$, a  
locally finitely presented Grothendieck category $\mathcal C$, which is usually  
called the functor category of ${\rm Mod}$-$R$, and a
fully faithful embedding
$$F:{\rm Mod-}R \rightarrow \mathcal{C}$$
satisfying the following properties:

\bigskip

\begin{enumerate}
\item The functor $F$ has a right adjoint functor $G:\mathcal{C} \rightarrow {\rm Mod}$-$R$.

\bigskip

\item A short exact sequence
$$\Sigma\equiv 0 \rightarrow X\rightarrow Y\rightarrow Z\rightarrow 0$$
in ${\rm Mod}$-$R$ is pure if and only if the induced sequence $F(\Sigma)$
is exact (and pure) in $\mathcal C$.

\bigskip

\item $F$ preserves finitely generated objects, i.e., the image of any finitely generated object in ${\rm Mod}$-$R$ is a finitely generated object in $\mathcal C$.

\bigskip

\item $F$ identifies ${\rm Mod}$-$R$ with the full subcategory of $\mathcal C$ consisting of the all
 FP-injective objects in $\mathcal C$ where an object  $C\in\mathcal C$ is FP-injective if ${\rm Ext}_{\mathcal C}^1(X,C)=0$ for every finitely presented object
 $X\in \mathcal C$.

\bigskip
 
\item A module $M\in{\rm Mod}$-$R$ is pure-injective if and only if $F(M)$ is an injective object of $\mathcal C$.

\bigskip

\item Every module $M\in{\rm Mod}$-$R$ admits a pure embedding in a pure-injective object $PE(M)\in{\rm Mod}$-$R$ such that the image of this embedding under
$F$ is the injective envelope of $F(M)$ in $\mathcal C$. The pure-injective object $PE(M)$ is called the
{\em pure-injective envelope} of $M$.
\end{enumerate}

\bigskip

\noindent We easily deduce the following key result from the above properties:

\begin{proposition}\label{pure}
Let $M$ be a right $R$-module. Then $M$ is invariant under automorphisms (resp., endomorphisms) of its pure-injective envelope 
in Mod-$R$ if and only if $F(M)$ is invariant under automorphisms (resp., endomorphisms) of its injective envelope in $\mathcal C$.
\end{proposition}

\noindent This proposition allows as to extend the result obtained for modules which are invariant under automorphisms of their injective envelope to this new setting. 

On the other hand, let us recall that a module object $Q$ in a Grothendieck category is called quasi-injective if every morphism $f:C\rightarrow Q$, where $C$ is a subobject of $Q$, extends to an endomorphism of $Q$. And this property is equivalent to the claim that $Q$ is invariant under endomorphisms of its injective envelope $E(Q)$. Therefore, it seems natural to state the following definition:

\begin{mydef}
A module $M$ in Mod-$R$ is called quasi pure-injective if it is invariant under endomorphisms of its pure-injective envelope, equivalently, if $F(M)$ is a quasi-injective object in the associated functor category $\mathcal C$.
\end{mydef}

\noindent It is well-known that any object $C$ in a Grothendieck category admits a minimal embedding $u:C\rightarrow Q$ in a quasi-injective object $Q$, which is called its quasi-injective envelope. In particular, this shows that, for any module $M$, the object $F(M)$ has a quasi-injective envelope $u:F(M)\rightarrow Q$ in the functor category $\mathcal C$ of Mod-$R$. By construction, $F(M)$ is an FP-injective object of $\mathcal C$, but  we cannot see any reason why this object $Q$ must be also FP-injective and thus, belong to the image of the functor $F$. As a consequence, it seems that in general, modules in Mod-$R$ do not need to have a quasi pure-injective envelope. 

Surprisingly, it seems that this natural problem has never been considered and clarified in the literature. Our next example tries to shed some light to the possible consequences that this lack of quasi pure-injective envelopes may have in the characterization of quasi pure-injective modules.

\begin{example}
Let $R$ be a commutative PID. Then an $R$-module $M$ is flat if and only if it is torsion-free and therefore, the only possible pure ideals of $R$ are $0$ and $R$ itself. This means that any homomorphism $f:N\rightarrow R$ from a pure ideal $N$ of $R$ to $R$ trivially extends to an endomorphism of $R$. However, if $R$ would always be quasi pure-injective, then we would get that 
$$R\cong \End_R(R) \cong \End_{\mathcal C}(F(R))\cong \End_{\mathcal C}(E(F(R)))$$
where $E(F(R))$ is the injective envelope of $F(R)$ in $\mathcal C$. And, as the endomorphism ring of an injective object in a Grothendieck category is always left pure-injective, this would mean that any commutative PID is pure-injective. But this is not the case, as for instance the case of $\mathbb Z$ shows.   
\end{example}

\noindent We finish this section by showing different applications of Proposition \ref{pure} to the study of modules $M$ which are invariant under automorphisms of their pure-injective envelope. Proofs are based on applying the corresponding results to the image $F(M)$ of the module in the functor category.  

\bigskip 

\noindent As an analogue of Theorem \ref{quasi}, we have the following in this new setting.

\begin{theorem} \cite{GKS}
Let $M$ be a right $R$-module. If $\End(M)$ has no homomorphic images isomorphic to $\mathbb{F}_2$, then $M$ is invariant under automorphisms of its pure-injective envelope if and only it is quasi pure-injective. 

In particular, a module $M$ over an algebra $A$ over a field with more than two elements is invariant under automorphisms of its pure-injective envelope if and only if it is quasi pure-injective.
\end{theorem}

\noindent The next theorem is an analogue of Theorems \ref{decomp} and \ref{struct}.

\begin{theorem} \cite{GKS}
Let $M$ be a module which is invariant under automorphisms of its pure-injective envelope. Then $M=A\oplus B$, where $A$ is quasi pure-injective and $B$ is square-free. 

Moreover, $\End(M)/J(\End(M))\cong R_1\times R_2$, where $R_1=\End(A)/J(\End(A))$ is von Neumann regular and right self-injective and $R_2=\End(B)/J(End(B))$ is abelian regular.
\end{theorem}

\noindent Since it is known that a quasi pure-injective module satisfies the full exchange property \cite[Theorem 11]{ZZ1} and any square-free module with the finite exchange property satisfies the full exchange property, we have the following as a consequence of the above theorem.

\begin{theorem} \cite{GKS}
Let $M$ be a module which is invariant under automorphisms of its pure-injective envelope. Then $M$ satisfies the full exchange property.
\end{theorem}
 
 \noindent Also, we have the following
 
\begin{theorem} \cite{GKS}
Let $M$ be a module which is invariant under automorphisms of its pure-injective envelope. Then $M$ is a clean module.  
\end{theorem}

\bigskip

\bigskip


\bigskip

\bigskip

\begin{thebibliography}{99}

\bigskip

\bigskip

\bibitem{AA1} 
G. Abrams and G. Aranda Pino, {\it The Leavitt path algebra of a graph}, J. Algebra 293, 2 (2005), 319-334.

\bibitem{AA2} 
G. Abrams and G. Aranda Pino, {\it The Leavitt path algebras of arbitrary graphs}, Houston J. Math. 34, 2 (2008), 423-442.

\bibitem{AAJZ}
A. Alahmadi, H. Alsulami, S. K. Jain and E. Zelmanov, {\it Structure of Leavitt path algebras of polynomial growth}, Proc. Natl. Acad. Sci. 110, 38 (2013), 15222-15224. 

\bibitem{AEJ} 
A. Alahmadi, N. Er and S. K. Jain, {\it Modules which are invariant under monomorphisms of their injective hulls}, J. Aust. Math. Soc 79 (2005), 349-360.

\bibitem{ARS}
G. Aranda Pino, K. M. Rangaswamy and M. Siles Molina, {\it Endomorphism rings of Leavitt path algebras}, preprint.

\bibitem{CH1} 
J. Clark and D. V. Huynh, {\it Simple rings with injectivity conditions on one-sided ideals}, Bull. Austral. Math. Soc. 76 (2007), 315-320.

\bibitem{CJ}
P. Crawley and B. J\'{o}nnson, {\it Refinements for infinite direct decompositions of algebraic systems}, Pacific J. Math. 14, 3 (1964), 755-1127.

\bibitem{CB} W. Crawley-Boevey, {\it Locally finitely presented additive categories}, Comm. Algebra 22 (1994), 1641-1674.

\bibitem{hai} H. Q. Dinh, {\it A note on pseudo-injective modules}, Comm. Algebra 33 (2005), 361-369.

\bibitem{DF}
S. E. Dickson and K. R. Fuller, {\it Algebras for which every indecomposable right module is invariant in its injective envelope}, Pacific J. Math. 31, 3 (1969), 655-658.

\bibitem{ESS}
N. Er, S. Singh and A. K. Srivastava, {\it Rings and modules which are stable under automorphisms of their injective hulls}, J. Algebra 379 (2013), 223-229.

\bibitem{FU}
C. Faith and Y. Utumi, {\it Quasi-injective modules and their endomorphism rings}, Arch. Math. 15 (1964), 166-174.

\bibitem{Fuchs}
L. Fuchs, {\it On quasi-injective modules}, Ann. Scuola Norm. Sup. Pisa 23 (1969), 541-546.

\bibitem{GH} P. A. Guil Asensio and I. Herzog, {\it Left cotorsion rings}, Bull. London Math.
Soc. 36 (2004), 303-309.

\bibitem{GKS}
P. A. Guil Asensio, D. Keskin T\"ut\"unc\"u and A. K. Srivastava, {\it Modules invariant under automorphisms of their covers and envelopes}, Israel J. Math., to appear.

\bibitem{GS1}
P. A. Guil Asensio and A. K. Srivastava, {\it Additive unit representations in endomorphism rings and an extension of a result of Dickson and Fuller}, Ring Theory and Its Applications, Contemp. Math., Amer. Math. Soc., 609 (2014), 117-121.

\bibitem{GS2}
P. A. Guil Asensio and A. K. Srivastava,  {\it Automorphism-invariant modules satisfy the exchange property}, J. Algebra 388 (2013), 101-106.

\bibitem{H}
P. Hill, {\it Endomorphism ring generated by units}, Trans. Amer. Math. Soc. 141 (1969), 99-105.

\bibitem{JS} S. K. Jain and S. Singh, {\it On pseudo-injective modules and self pseudo injective rings}, J. Math. Sci. 2 (1967), 23-31.

\bibitem{JW}
R. E. Johnson and E. T. Wong, {\it Quasi-injective modules and irreducible rings}, J. London Math. Soc. 36 (1961), 260-268.


\bibitem{KS1}
D. Khurana and A. K. Srivastava, {\it Right self-injective rings in which
each element is sum of two units}, J. Algebra Appl. 6, 2 (2007), 281-286.

\bibitem{KS2}
D. Khurana and A. K. Srivastava, {\it Unit sum numbers of right
self-injective rings}, Bull. Australian Math. Soc. 75, 3 (2007), 355-360.

\bibitem{LZ2}
T. K. Lee and Y. Zhou, {\it On lifting of idempotents and semiregular endomorphism rings}, Colloq. Math. 125, 1 (2011), 99-113.

\bibitem{LZ}
T. K. Lee and Y. Zhou, {\it Modules which are invariant under automorphisms of their injective hulls}, J. Algebra Appl. 12 (2), (2013).

\bibitem{MM} S. H. Mohamed and B. J. M\"{u}ller, {\it Continuous and Discrete Modules}, Cambridge University Press, 1990.

\bibitem{Nichol}
W. K. Nicholson, {\it Lifting idempotents and exchange rings}, Trans. Amer. Math. Soc. 229, (1977), 269-278.

\bibitem{O}
B. L. Osofsky, {\it Endomorphism rings of quasi-injective modules}, Canad. J. Math. 20 (1968), 895-903.

\bibitem{Si} D. Simson, {\it On pure semi-simple Grothendieck 
categories I}, Fund. Math. 100 (1978),
211-222.

\bibitem{St} R.  W.  Stringall,  {\it Endomorphism  rings  of  abelian  groups  generated  by  automorphism  groups},
Acta. Math. 18 (1967), 401-404.

\bibitem{SS1}
S. Singh and A. K. Srivastava, {\it Dual automorphism-invariant modules}, J. Algebra 371 (2012), 262-275.

\bibitem{SS2}
S. Singh and A. K. Srivastava, {\it Rings of invariant module type and automorphism-invariant modules}, Ring Theory and Its Applications, Contemp. Math., Amer. Math. Soc., 609 (2014), 299-311.

\bibitem{Teply}
M. L. Teply, {\it Pseudo-injective modules which are not quasi-injective}, Proc. Amer. Math. Soc. 49, 2 (1975), 305-310.

\bibitem{Warfield1}
R. B. Warfield, Jr., {\it Decompositions of injective modules}, Pacific J. Math. 31 (1969), 263-276.

\bibitem{Warfield2}
R. B. Warfield, Jr., {\it Exchange rings and decompositions of modules}, Math. Ann. 199 (1972), 31-36. 

\bibitem{ZZ}
B. Huisgen-Zimmermann and W. Zimmermann, {\it Algebraically compact rings and modules}, Math. Z. 161 (1978), 81-93.

\bibitem{ZZ1} B. Huisgen-Zimmermann and W. Zimmermann, {\it 
Classes of modules with the exchange property},  
J. Algebra 88, 2 (1984), 416-434.

\end{thebibliography}
\end{document}